 \documentclass[final,3p,times,twocolumn]{elsarticle}
\usepackage{amssymb}
\usepackage{amsthm}
\usepackage{amsmath,amsfonts}

\newtheorem{theorem}{Theorem}

\usepackage{algorithmic}
\usepackage{algorithm}

\usepackage{xcolor}

\biboptions{sort&compress}

\usepackage{graphicx}
\usepackage[colorlinks=true, allcolors=blue]{hyperref}
\DeclareGraphicsExtensions{.eps,.ps,.jpg,.bmp}
\usepackage{subfigure}

\usepackage{lineno}
\modulolinenumbers[5]

\journal{}

\begin{document}

\begin{frontmatter}


 \title{Solving a class of stochastic optimal control problems by physics-informed neural networks\tnoteref{label1}}
 \tnotetext[label1]{This research was partially supported by the National Natural Science Foundation of China (12272297).}
 
  \author[inst1]{Zhe Jiao\fnref{equal}}
  \ead{zjiao@nwpu.edu.cn}
  
 \author[inst1,inst2]{Wantao Jia\corref{cor1}\fnref{equal}}
 \ead{jiawantao@nwpu.edu.cn}
 \cortext[cor1]{Corresponding author}
 
\author[inst3]{Weiqiu Zhu}
\ead{wqzhu@zju.edu.cn}
 
 \affiliation[inst1]{organization={School of Mathematics and Statistics},
             addressline={Northwestern Polytechnical University},
             city={Xi'an},
             postcode={710129},
             country={China}}
             
\affiliation[inst2]{organization={MOE Key Laboratory for Complexity Science in Aerospace},%Department and Organization
             addressline={Northwestern Polytechnical University},
             city={Xi'an},
             postcode={710129},
             country={China}}
             
\affiliation[inst3]{organization={State Key Laboratory of Fluid Power and Mechatronic Systems, Department of Mechanics},%Department and Organization
             addressline={Zhejiang University},
             city={Hangzhou},
             postcode={310027},
             country={China}}

\fntext[equal]{Equal contribution.}
             
%\title{Title of Your Manuscript}
%
%%% use optional labels to link authors explicitly to addresses:
%%% \author[label1,label2]{}
%%% \affiliation[label1]{organization={},
%%%             addressline={},
%%%             city={},
%%%             postcode={},
%%%             state={},
%%%             country={}}
%%%
%%% \affiliation[label2]{organization={},
%%%             addressline={},
%%%             city={},
%%%             postcode={},
%%%             state={},
%%%             country={}}
%
%\author[inst1]{Author One}
%
%\affiliation[inst1]{organization={Department One},%Department and Organization
%            addressline={Address One}, 
%            city={City One},
%            postcode={00000}, 
%            state={State One},
%            country={Country One}}
%
%\author[inst2]{Author Two}
%\author[inst1,inst2]{Author Three}
%
%\affiliation[inst2]{organization={Department Two},%Department and Organization
%            addressline={Address Two}, 
%            city={City Two},
%            postcode={22222}, 
%            state={State Two},
%            country={Country Two}}

\begin{abstract}
The aim of this work is to develop a deep learning method for solving high-dimensional stochastic control problems based on the Hamilton--Jacobi--Bellman (HJB) equation and physics-informed learning. 
Our approach is to parameterize the feedback control and the value function using a decoupled neural network with multiple outputs. 
We train this network by using a loss function with penalty terms that enforce the HJB equation along the sampled trajectories generated by the controlled system.  
More significantly, numerical results on various applications are carried out to demonstrate that the proposed approach is efficient and applicable. 
\end{abstract}




\end{frontmatter}

%\linenumbers

%%%%%%%%%%%%%%%%%%%%%%%%%%%
%%%%%%%%%%%%%%%%%%%%%%%%%%%
\section{Introduction}
\label{sec:introduction}
The range of stochastic optimal control (SOC) problems covers a variety of scientific branches such as finance~\cite{pham2009continuous}, molecular dynamics \cite{gao2023transition}, neuroscience \cite{Todorov2004} and robotics \cite{manipulation}.
To address SOC problems, there are two prominent frameworks: Pontryagin's maximum principle~(MP)~\cite{Pontrygin1987} and Bellman's dynamic programming~(DP)~\cite{Bellman1958}.  
Drawing on these frameworks, many numerical methods have been developed for tackling SOC problems (cf. \cite{kushner1990numerical,jin2022survey} and references therein). 

%\begin{table*}[htbp]
%%\vskip 0.15in
%\begin{center}
%\begin{tabular}{cclccc}
%\hline
%Papers & Theory & Reformulation of SOC & Networks & $\star$ \\
%\hline
%\cite{HJE2018, GPW2022} & DP & backward SDE & FC & no \\
%\cite{JPPZ2022} & MP & Hamiltonian system & FC & no \\
%\cite{HPBL2021, HPBL2022} & DP &State-action value function & FC & no\\
%\cite{PWCRT2020, Wang2019}& DP & forward-backward SDEs & LSTM & yes\\
%This work & DP & pathwise HJB operator & PINN & no \\
%\hline
%\end{tabular}
%\end{center}
%\caption{Comparison of this paper to some related works. Here, $\star$ denotes whether the optimal controls have an explicit representation.}
%\label{comparison}
%%\vskip 0.1in
%\end{table*}

However, these traditional numerical methods are not applicable when the state dimension is large \cite{ET2018}. In recent years, there has seen significant progress in leveraging deep learning (DL) to solve the high-dimensional SOC problems \cite{GKM2015,HE2016,Wang2019,li2024neural,cai2024soc}.
Broadly speaking, the deep neural network-based methods for SOC can be divided into two distinct categories. In the first category, it is concerned with the DL-based approach to solve the extended Hamiltonian system, which is derived from stochastic MP (cf. \cite{FZ2020,CL2022,JPPZ2022}). For the study of the second category, \cite{HPBL2021, HPBL2022} reformulate the SOC problem as Markov decision process based on DP, which is solved by some DL-based algorithms. Another direction of this category is to solve the SOC problem from the view of DP via HJB equation \cite{PWCRT2020,nusken2021solving,hua2024simulation}. We need to point out that in these papers Feynman--Kac formula is the basis to probabilistically represent the solution to HJB equation so that the author can utilizes neural networks to obtain the optimal policy.     

Motivated by previous research, we aim to solve the SOC problem with physics-informed learning~\cite{RPK2019,karniadakis2021physics}. The main issue we encounter in our approach is to construct a physics-informed neural network (PINN) for solving HJB equation, which is a semilinear parabolic partial differential equation (PDE) with a terminal value condition. Since the HJB equation is defined on the whole space, without boundary condition, PINN cannot be directly used to compute the value function by solving the HJB equation. Thanks to the stochastic verification theorem (see Theorem~\ref{verification} in Section~\ref{sec:verification}), we can simulate the value function along the trajectories of the controlled system, not on the whole space, by neural network. This is the key idea of our approach.   

Our main contribution is twofold: 
(i) In contrast to~\cite{li2024neural}, we use the controlled SDE to conduct sampling on relevant states during PINN training;
(ii) We propose a simulation-free algorithm for SOC by physics-informed learning, which means it dose not require numerical solutions of the control problem.

The remaining part of this paper is organized as follows. 
In Section \ref{sec:soc}, we briefly introduce the preliminaries about the SOC problem, and the verification theorem that is the basis to construct our DL based solver. This solver called DeepHJB is proposed in Section~\ref{sec:deepHJB}.
Numerical examples in Section \ref{sec:experiments} illustrate our proposed solver to solve some SOC problems. 
Section \ref{sec:conclusion} provides some conclusions.
%%%%%%%%%%%%%%%%%%%%%%%%%%%
%%%%%%%%%%%%%%%%%%%%%%%%%%%
\section{Stochastic optimal control}
\label{sec:soc}

\subsection{Problem setup}
\label{sec:problem}
Let $T>t\geqslant 0$ and $\textbf{W}: [t, T]\times\Omega\rightarrow \mathbb{R}^{d}$ be a $d$-dimensional standard $\mathbb{F}$-Brownian motions on a filtered probability space $(\Omega, \mathcal{F}, \mathbb{F}, \mathbb{P})$ where $\mathbb{F}= \{\mathcal{F}_s\}_{t\leqslant s \leqslant T}$ is the natural filtration generated by $\textbf{W}(s)$. The quadruple $(\Omega, \mathcal{F}, \mathbb{F}, \mathbb{P})$ also satisfies the usual hypotheses (see Chapter 1.4 in~\cite{chung2013}). $\mathbb{E}[\cdot]$ stands for expectation with respect to the probability measure $\mathbb{P}$.

We consider the controlled stochastic differential equation (SDE) as follows
\begin{equation}\label{SDE}
   \mathrm{d}\textbf{x}_s = b(s, \textbf{x}_s, \textbf{u}(s))\mathrm{d}s + \sigma(s, \textbf{x}_s, \textbf{u}(s))\mathrm{d}\textbf{W}(s)
\end{equation}
with $s\in[t_0, T]$ and the initial data $\textbf{x}_{t_0} = x\in\mathbb{R}^n$. 
Here, $\textbf{x}_s \in\mathbb{R}^n$ is the state process, $\textbf{u}(s)\in\mathbb{R}^m$ is a control process valued in a given subset $U$ of $\mathbb{R}^{m}$.
The cost functional is given by
\begin{equation}\label{cost}
{\small
\begin{aligned}
   J(t_0, x; \textbf{u}(t)) = \mathbb{E}\left[ \int_{t_0}^{T}\phi(s, \textbf{x}_s, \textbf{u}(s))\mathrm{d}s + \psi(\textbf{x}_T) | \textbf{x}_{t} = x \right]
\end{aligned}
}
\end{equation}
with the functions $\phi: [t_0, T]\times\mathbb{R}^{n} \times \mathbb{R}^m \rightarrow \mathbb{R}$ and $\psi: \mathbb{R}^{n}\rightarrow \mathbb{R}$.
The goal of our SOC problem is to look for an admissible control (if exists) that minimizes (\ref{cost}) over $\mathcal{U}$ which is the set of all admissible controls defined by
\[
	\mathcal{U}:=\left\{u: [t_0, T] \times \Omega \rightarrow U| u(s) \in L^2_{\mathbb{F}}(t_0, T; \mathbb{R}^{m})\right\}
\]
in which $L^2_{\mathbb{F}}(t, T; \mathbb{R}^{m})$ consists of all $\mathbb{F}$-adapted functions $u: [t_0, T]\times\Omega\rightarrow \mathbb{R}^{m}$ satisfying
$\mathbb{E}[ \int_{t_0}^{T}|u|^2\mathrm{d}s]<\infty$.

In this paper, we focus on the SOC problem under the following conditions.
\begin{itemize}
\item
The drift term $b\in\mathbb{R}^n$ and the diffusion term $\sigma\in\mathbb{R}^{n\times(1+d)}$ in (\ref{SDE}) have the following linear forms in control
\begin{equation*}
b(s, \textbf{x}_s, \textbf{u}(s)) = A(s, \textbf{x}_s) + B(s, \textbf{x}_s) \textbf{u}(s),
\end{equation*}
and
\begin{equation*}
\sigma(s, \textbf{x}_s, \textbf{u}(s)) = [\lambda B(s, \textbf{x}_s) \textbf{u}(s), C(s, \textbf{x}_s)]
\end{equation*}
with $\lambda \geqslant 0$, $A\in\mathbb{R}^{n}$, $B\in\mathbb{R}^{n\times m}$ and $C \in \mathbb{R}^{n\times d}$.
\item
The random term $\textbf{W}_s = [w^{(1)}_s, w^{(2)}_s]\in\mathbb{R}^{(1+d)}$ in which $w^{(1)}_s\in\mathbb{R}^{1}$ and $w^{(2)}_s\in\mathbb{R}^{d}$ are mutually independent Brownian motions. 
\item
The running cost in (\ref{cost}) is quadratic, that is, 
\[
	\phi(s, \textbf{x}_s, \textbf{u}(s)) = \textbf{x}_s^\top F\textbf{x}_s+ \frac{1}{2}\textbf{u}(s)^{\top}D\textbf{u}(s)
\]
with the coefficients $F\in\mathbb{R}^{n\times n}$ and $D\in\mathbb{R}^{m\times m}$.
\item
The terminal cost is linear 
\[
	\psi(x) = \gamma \cdot x
\]
or quadratic
\[
	\psi(x) = (x-\textbf{x}_T)^\top F_T (x-\textbf{x}_T)
\]
with the coefficients $\gamma\in\mathbb{R}^n$ and $F_T\in\mathbb{R}^{n\times n}$.
\end{itemize}
Under suitable assumptions (see Chapter 1 in \cite{YZ1991}),  for any $\textbf{u}(s)\in\mathcal{U}$ equation~\eqref{SDE} has a unique solution $\textbf{x}_s$ and the cost function~\eqref{cost} is well-defined. We call $(\textbf{x}_s, \textbf{u}(s))$ an admissible pair.	
Any $\textbf{u}^\ast(s)$ is called an optimal control if it satisfies
\begin{equation*}
    \textbf{u}^\ast(s) := \mathop{\arg\min}\limits_{\mathbf{u}(s)\in\mathcal{U}}J(t, x; \textbf{u}(s)).
\end{equation*}
The corresponding state process $\textbf{x}^{\ast}_s$ is called an optimal trajectory and the state-control pair $(\textbf{x}^{\ast}_s, \textbf{u}^\ast(s))$ called an optimal pair.

\subsection{Verification theorem}
\label{sec:verification}
We define the value function as
\begin{equation*} \label{valuefunction}
	 q(t, x):= J(t, x; \textbf{u}^\ast(s)) = \min\limits_{\mathbf{u}(s)\in\mathcal{U}} J(t, x; \textbf{u}(s)).
\end{equation*}
The following theorem shows the evolution of the value function along the optimal trajectory and at the the final time, which is deduced from the stochastic verification theorem (see Theorem 5.1 in Chapter 5.5 of \cite{YZ1991}). The detailed proof is given in~\ref{sec:thmproof}.

\begin{figure*}[!htbp]
\centering
\includegraphics[width=0.8\textwidth]{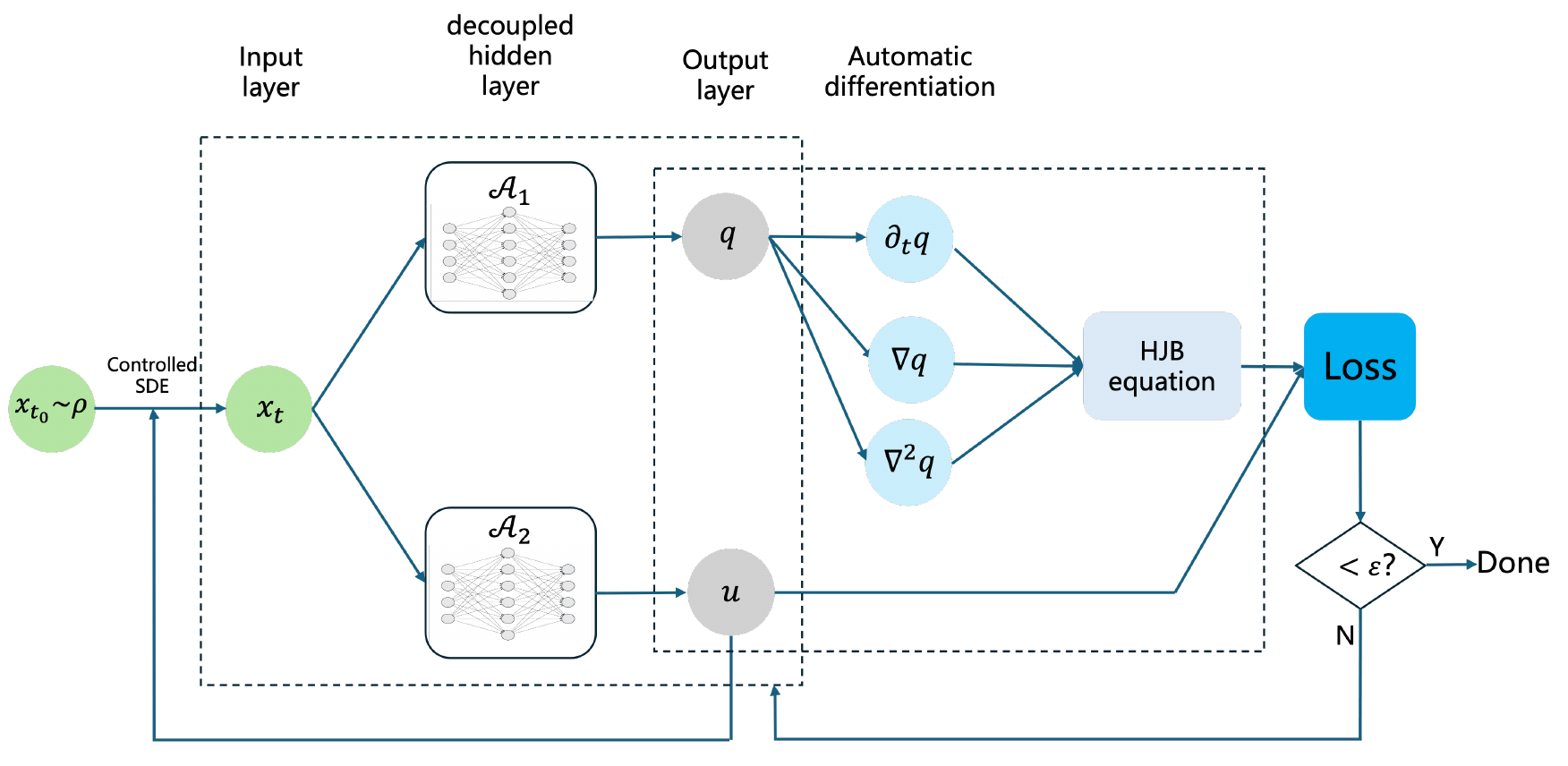}
\caption{A decoupled neural network structure for solving HJB equation~\eqref{discreteHJB} with multiple outputs. This neural network is a type of PINNs which integrate the information from PDEs~\eqref{discreteHJB} into the loss function of a neural network using automatic differentiation. The architecture of the decoupled hidden layers are denoted by $\mathcal{A}_i$, $i=1,2$, which will be given in~\ref{sec:configuration}.}
\label{decoupled_nn}
\end{figure*}
\begin{theorem}\label{verification}
An admissible pair $(\mathbf{x}_t, \mathbf{u}(t))$, where the feedback control $\mathbf{u}(t)$ is given by 
\begin{equation} \label{optimalcontrol}
	\mathbf{u}(t) = -\tilde{D}(t, \mathbf{x}_t)^{-1}B(t, \mathbf{x}_t)^{\top } \nabla q(t, \mathbf{x}_t)
\end{equation}
with 
\[
	\tilde{D}(t, \mathbf{x}_t)=D+ \lambda^2 B(t, \mathbf{x}_t)^{\top}\nabla^2 q(t, \mathbf{x}_t) B(t, \mathbf{x}_t),
\]
is optimal if and only if the following {\rm{HJB}} equation holds
\begin{equation} \label{pathwiseValue}
\begin{aligned}
    	 -\partial_t q(t, \mathbf{x}_t) = &H\left(t, \mathbf{x}_t, \mathbf{u}(t), \nabla q(t, \mathbf{x}_t), \nabla^2 q(t, \mathbf{x}_t)\right) \\
	 =&\frac{1}{2}\mathrm{tr}\left[C(t, \mathbf{x}_t)^{ \top }\nabla^2q C(t, \mathbf{x}_t) \right]\\
	 & -\frac{1}{2}(\nabla q)^{\top} \left[B(t, \mathbf{x}_t)\tilde{D}^{-1}B(t, \mathbf{x}_t)^{ \top }\right]\nabla q\\
	 &+ A(t, \mathbf{x}_t)\cdot \nabla q + \mathbf{x}_t^\top F\mathbf{x}_t
\end{aligned}
\end{equation}
for any $t\in[0, T)$, and $q(T, \mathbf{x}_{T}) = \psi(\mathbf{x}_T)$.
\end{theorem} 
Here, $\partial_t q$ means the first-order derivative of $q$ with respect to $t$, $\nabla q$ and $\nabla^2 q$ respectively denote the gradient and the Hessian of $q(t, x)$ with respect to $x$, and $\mathrm{tr}$ is the abbreviation of the trace operator.

\section{Deep learning approach}
\label{sec:deepHJB}
In this section we propose our approach to seek an optimal pair that minimizes the cost functional~\eqref{cost} subject to~\eqref{SDE} for initial data sampled from a probability distribution in $\mathbb{R}^n$ with a density denoted by $\rho$.

We select a partition of the time interval $[0, T]$: 
{\small
\[
	0=t_0 < t_1 < \cdots< t_n<\cdots <t_N =T,
\]}
and denote by $\triangle t_n = t_{n+1}-t_{n}$ the $(i+1)$th interval of the grid and $\triangle W_n = W_{t_{n+1}} - W_{t_{n}}$ the $(i+1)$-th increment of the Brownian motion.
Once the control $u(t)$ is computed, the Euler--Maruyama scheme (cf. \cite{KP}) of (\ref{SDE}) gives
\begin{equation} \label{scheme}
%{\footnotesize
\begin{aligned}
    x_{t_{n+1}} - x_{t_{n}} = b(t_n, x_{t_{n}}, u(t_{n}))\triangle t_n + \sigma(t_n, x_{t_{n}}, u(t_{n}))\triangle W_n
\end{aligned}
%}
\end{equation}
with the initial data $x_{t_{0}}= x\sim\rho$.
Using the numerical scheme~\eqref{scheme}, the path $\{(t_n, x_{t_{n}})\}_{0\leqslant n \leqslant N}$ can be easily generated. If the value function $q(t, x_t)$ is known and the discretization of the admissible pair $\{(x_{t_{n}}, u(t_{n}))\}_{0\leqslant n \leqslant N}$ satisfies
    \begin{equation}\label{discreteHJB}
    %{\footnotesize
    \left\{
                \begin{array}{ll}
      -\partial_t q(t_n, x_{t_n})= H\left(t_n, x_{t_n}, u(t_n), \nabla q(t_n, x_{t_n}), \nabla^2 q(t_n, x_{t_n})\right) ,\\
        q\left(t_N, x_{t_{N}}\right) = \psi\left(x_{t_{N}}\right),
                \end{array}
        \right.
        %}
    \end{equation}
from \eqref{pathwiseValue} in Theorem \ref{verification} we know $\{(x_{t_{n}}, u(t_{n}))\}_{0\leqslant n \leqslant N}$ is an optimal pair. 

Our approach parameterizes the functions $u$ and $q$ by a decoupled neural network (Figure~\ref{decoupled_nn}), which are given by
{\small
\begin{equation*}
u^{\textrm{NN}}(t_{n}, x_{t_{n}}; \theta_u), \quad q^{\textrm{NN}}(t_{n}, x_{t_{n}}; \theta_q). 
\end{equation*}}
We denote by $\Theta=\{\theta_u, \theta_q\}$ the weights of the neural network, which is trained by minimizing the sum of the expected losses that arises from the following penalty terms.
\begin{itemize}
	\item The second-order HJB penalty terms are defined as
		\begin{equation*}
		{\small
		\begin{aligned}
			\mathcal{L}_1(\theta_u, \theta_q)=\left|\partial_t q^{\textrm{NN}} + H\left(t_n, x_{t_n}, u^{\textrm{NN}}, \nabla q^{\textrm{NN}}, \nabla^2 q^{\textrm{NN}}\right)\right|
		\end{aligned}
		}
		\end{equation*}
		and
		\begin{equation*}
		{\small
		\begin{aligned}
			\mathcal{L}_2(\theta_q) = &\left|q^{\textrm{NN}}\left(t_N, x_{t_{N}}; \theta_q\right) - \psi\left(x_{t_{N}}\right)\right|,
		\end{aligned}
		}
		\end{equation*}
%		\begin{equation*}
%		\begin{aligned}
%			\mathcal{L}_2(\theta_q, \theta_Q) = &\left|q^{\textrm{NN}}\left(t_N, x_{t_{N}}; \theta_q\right) - \psi\left(x_{t_{N}}\right)\right|\\
%			&+ \left|Q^{\textrm{NN}}\left(t_N, x_{t_{N}}; \theta_Q\right) -  \nabla\psi\left(x_{t_{N}}\right)\right|,
%		\end{aligned}
%		\end{equation*}
		where $\mathcal{L}_1$ and $\mathcal{L}_2$ are derived from~\eqref{discreteHJB}.
	\item Another penalty term is given as follows
	{\small
	\[
		\mathcal{L}_3(\theta_u, \theta_q) =\left|u^{\textrm{NN}} +  \tilde{D}^{-1}B^{\top }(t_n, x_{t_n}) \nabla q^{\textrm{NN}}\right|
	\]}
	with
	{\small
	\[
		\tilde{D}=D+ \lambda^2 B(t_n, x_{t_n})^{\top}\nabla^2 q^{\textrm{NN}}(t_n, x_{t_n}; \theta_q) B(t_n, x_{t_n}),
	\]}
	where $\mathcal{L}_3$ is from~\eqref{optimalcontrol}.
\end{itemize}
Now, we can define the \textit{physics-informed learning problem} 
\[
	\min_{\Theta}\mathbb{E}_{x\sim\rho}\left\{\alpha_1\mathcal{L}_1 + \alpha_2 \mathcal{L}_2 + \alpha_3 \mathcal{L}_3\right\}.
\]
The coefficients $\alpha_1>0$, $\alpha_2>0$ and $\alpha_3> 0$ are supposed to be fixed. 
%When we have no information on the optimal control, $\alpha_3$ is set to be $0$.    

Finally, we apply a SGD-type algorithm to optimize the parameter $\Theta$. 
The pseudo-code for implementing the above approach is given in Algorithm \ref{alg1}. 

\begin{algorithm}
\caption{DeepHJB solver}
\label{alg1}
\textbf{Input:} {\footnotesize the initial data $\{(t_0, x^{(i)}_{t_{0}})\}_{1\leqslant i\leqslant M}$, parameter $N$ of partition, parameters $\theta^{(0)}$ of networks, learning rate $\eta$, max-step $K$}
                  
 \hrulefill
 \\  
\textbf{For} $k = 0$ to $K-1$ \textbf{do}

\quad \textbf{For} $i=1$ to $M$ \textbf{do}

\quad\quad \textbf{For} $n=0$ to $N$ \textbf{do}

\quad\quad\quad {\footnotesize$q^{\textrm{NN}}(t_{n}, x^{(i)}_{t_{n}}; \theta_q^{(k)})$}

\quad\quad\quad {\footnotesize$u^{\textrm{NN}}(t_{n}, x^{(i)}_{t_{n}}; \theta_u^{(k)})$}

\quad\quad\quad \textbf{while} $n+1\leqslant N$ \textbf{do}
 
\quad\quad\quad\quad {\footnotesize$\triangle t_n = t_{n+1}-t_{n}$}

\quad\quad\quad\quad {\footnotesize$\triangle W^{(i)}_n = W^{(i)}_{t_{n+1}} - W^{(i)}_{t_{n}}$}

\quad\quad\quad\quad {\footnotesize $x^{(i)}_{t_{n+1}} = x^{(i)}_{t_{n}} + b(t_n, x^{(i)}_{t_{n}}, u^{\textrm{NN}})\triangle t_n + \sigma(t_n, x^{(i)}_{t_{n}}, u^{\textrm{NN}})\triangle W^{(i)}_n $}

%\quad\quad\quad\quad \quad\quad \quad\quad$+ \sigma(t_n, x^{(i)}_{t_{n}}, u^{\textrm{NN}})\triangle W^{(i)}_n$

\quad\quad\quad \textbf{end while}

\quad\quad \textbf{end for}

\quad  \textbf{end for}

\quad  {\footnotesize$\Theta^{(k)}  = ( \theta_u^{(k)}, \theta_q^{(k)})$}

\quad  {\footnotesize random set $B_k \subset \{1, 2, \cdots, M\}$}

\quad {\footnotesize$\mathrm{Loss} = \frac{1}{|B_k|}\sum\limits_{i\in B_k}\left\{\alpha_1\mathcal{L}^{(i)}_1(\Theta^{(k)}) + \alpha_2 \mathcal{L}^{(i)}_2(\Theta^{(k)}) + \alpha_3 \mathcal{L}^{(i)}_3(\Theta^{(k)})\right\}$}

\quad  {\footnotesize$\Theta^{(k+1)} = \Theta^{(k)} - \eta\nabla \mathrm{Loss}$}

\textbf{end for}

% \hrulefill
% \\
%\textbf{Outut:} $\{q(t_n, x_{t_{n}})\}_{0\leqslant n \leqslant N}$, $\{u(t_n)\}_{0\leqslant n \leqslant N}$ 
\end{algorithm}

%%%%%%%%%%%%%%%%%%%%%%%%%%%
%%%%%%%%%%%%%%%%%%%%%%%%%%%
\section{Numerical experiments}
\label{sec:experiments}
In this section, we apply the DeepHJB solver to some SOC problems. In the following subsections, we discuss the controlled Ornstein--Uhlenbeck (OU) dynamics and the controlled metastable dynamics, respectively. To evaluate the proposed solver, we introduce the following $L^2$ error as the performance metric
\[
	\mathbb{E}\left[\int_{0}^{T}|u^{\textrm{NN}} - u^\ast|^2(t, x_t)ds\right],
\] 
where $u^\ast$ is the baseline optimal control. The detailed configurations of these
experiments can be seen in~\ref{sec:configuration}.

\subsection{Ornstein--Uhlenbeck dynamics}
\label{sec:OUprocess}
We investigate the controlled system with
\begin{equation*}
\begin{split}
	A &=-I_{n\times n} + (\xi_{ij})_{1\leqslant i, j\leqslant n}, \\ B=C &= I_{n\times n} + (\xi_{ij})_{1\leqslant i, j\leqslant n}, \quad \lambda = 0,
\end{split}
\end{equation*}
where $\xi_{ij}\sim \mathcal{N}(0, 0.01)$ are sampled once at the beginning of the experiments.

\begin{figure}[!htbp]
\centering
\subfigure[Optimal control]{
\includegraphics[width=0.3\textwidth]{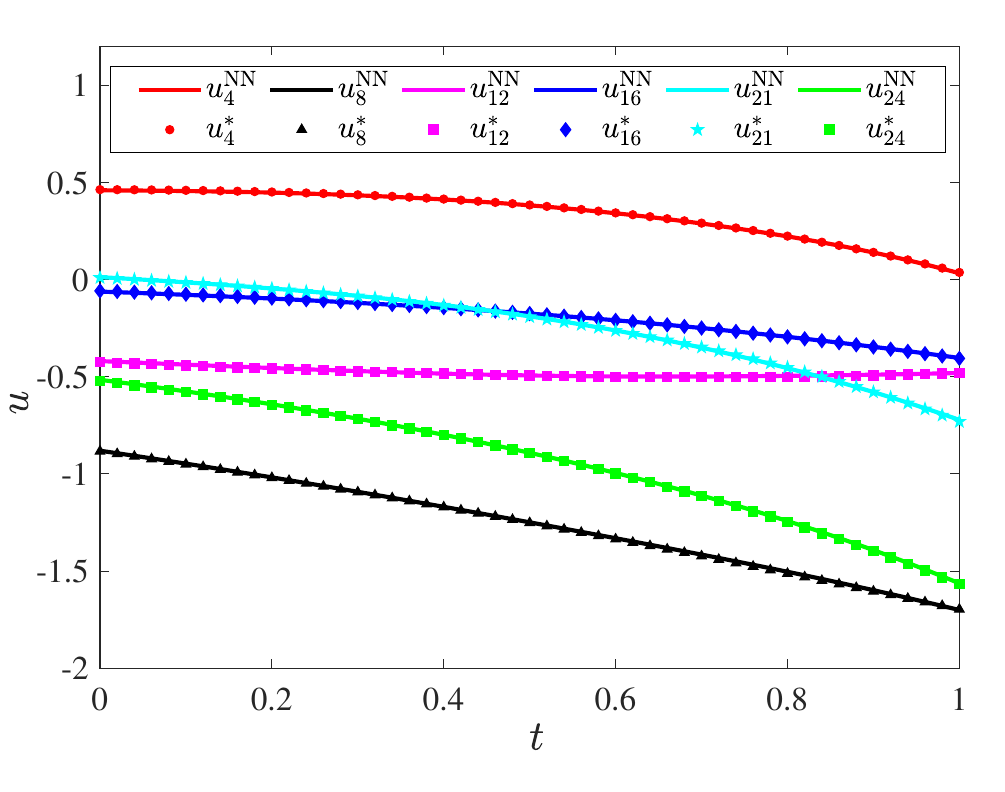}
}
\hspace{0in}
\subfigure[$L^2$ error]{
\includegraphics[width=0.3\textwidth]{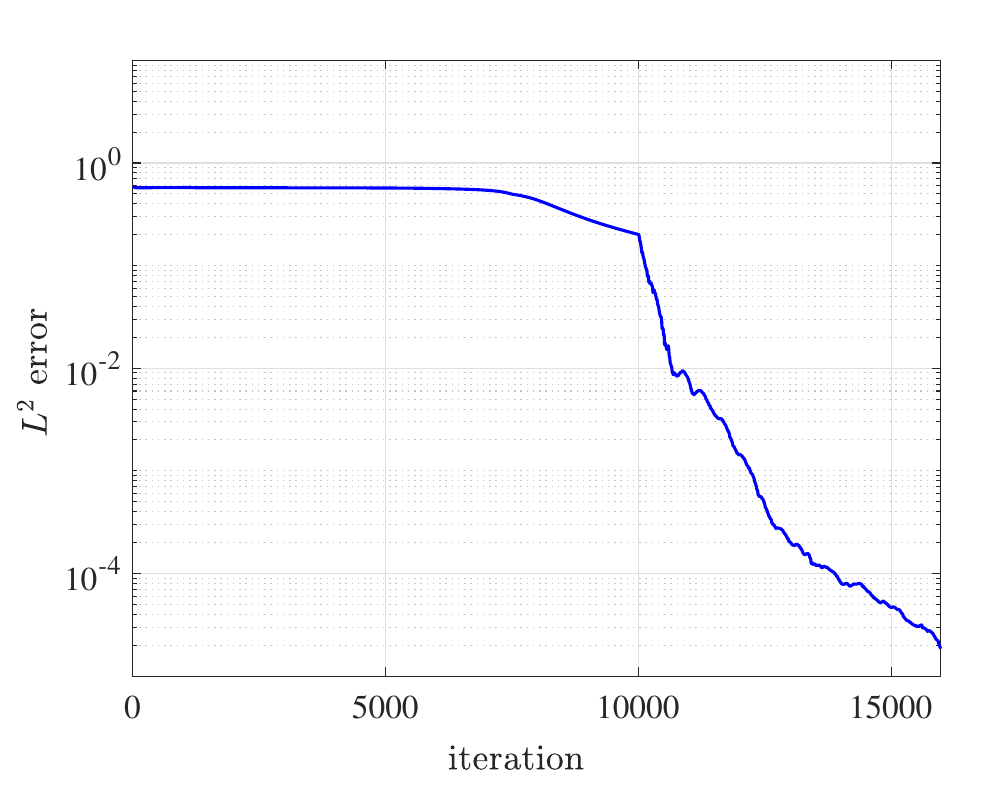}
}
\caption{Performance of the DeepHJB solver for $30$ dimensional controlled OU dynamics with linear terminal cost. (a) Comparison of $u^{\textrm{NN}}$ and $u^\ast$. Here, $7$ out of the $30$ components of the optimal control are pictured. (b) $L^2$ error with respect to the iteration step.}
\label{OU_1}
\end{figure}

For the SOC problem with linear terminal cost, we choose
\[
F=0,\quad D = I_{n\times n}, \quad \gamma = (1, \cdots, 1)^\top.
\]
In this situation, the optimal control can be given analytically by
\[
	u^\ast(t) = -B^{\top}e^{A^{\top}(T-t)}\gamma,
\]
which has been calculated in~\cite{nusken2021solving}.
We set the initial value to be zero and the terminal time $T=1.0$. In Figure~\ref{OU_1}, the subfigure (a) gives a visible comparison of the optimal control between $u^{\textrm{NN}}$ calculated by the DeepHJB solver and the baseline $u^\ast$, while the subfigure (b) shows the evolution of the error $L^2$ against the iteration step. It can be seen that the optimal control approximated by the proposed solver well coincides with the analytical one.

\begin{figure}[!htbp]
\centering
\subfigure[Optimal control]{
\includegraphics[width=0.3\textwidth]{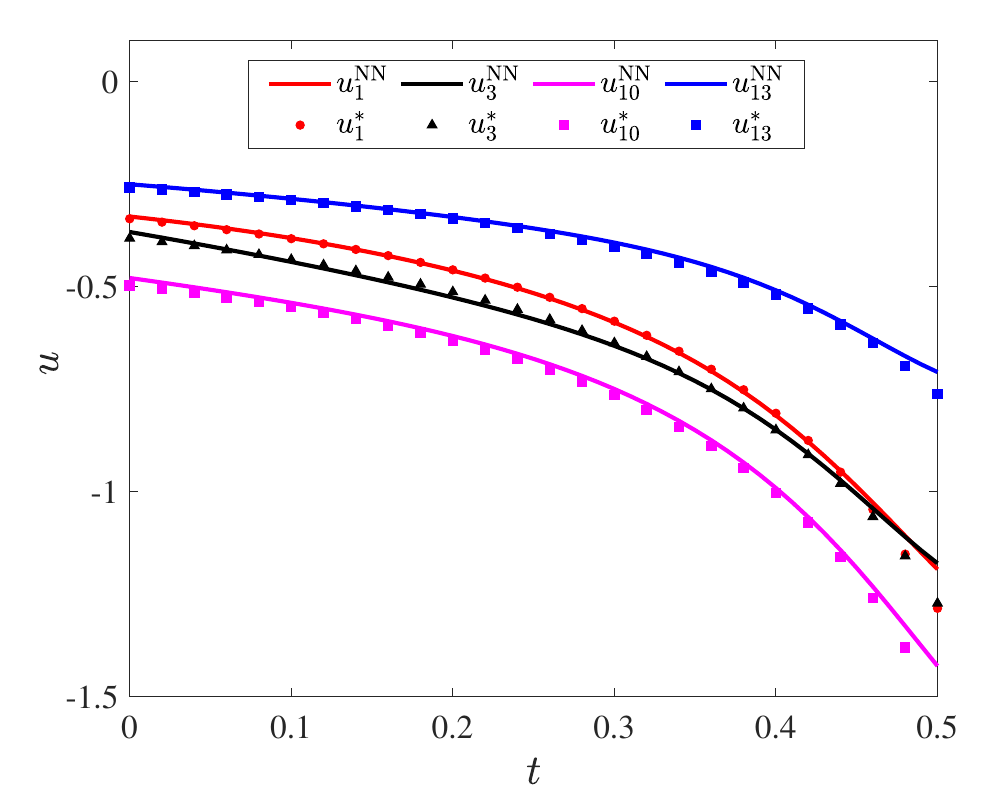}
}
\hspace{0in}
\subfigure[$L^2$ error]{
\includegraphics[width=0.3\textwidth]{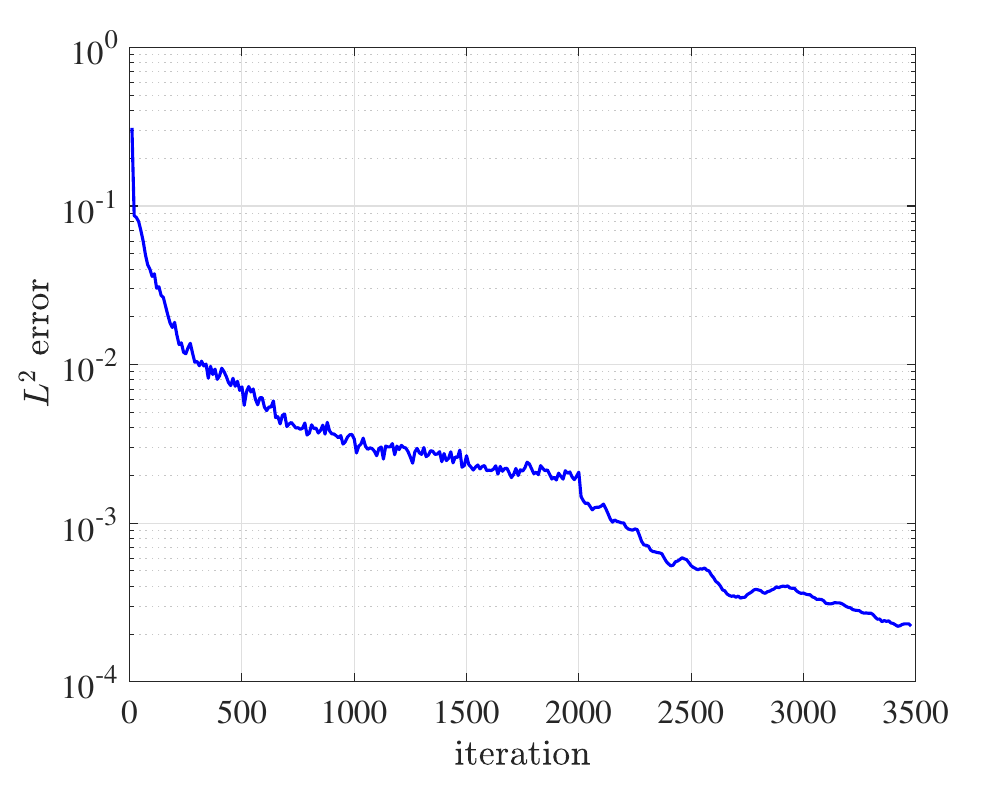}
}
\caption{Performance of the DeepHJB solver for $15$ dimensional controlled OU dynamics with a quadratic terminal cost. (a) Comparison of $u^{\textrm{NN}}$ and $u^\ast$ for varying time $t$ and at a fixed $x = (0.5, \cdots, 0.5)$. Here, four components of the optimal control are plotted. (b) $L^2$ error with respect to the iteration step.}
\label{OU_2}
\end{figure}

Regarding the case with a quadratic terminal cost, we choose
\begin{equation*}
\begin{split}
	F= \frac{1}{2}I_{n\times n}, \quad D = I_{n\times n}, \quad F_{T} = I_{n\times n}.
\end{split}
\end{equation*} 
This type of problems has an analytic optimal control
\begin{equation*}
\mathbf{u}^{\ast}(t, x) = -2B^{\top}P_{t}x
\end{equation*}
in which $P_{t}$ fulfills the Riccati equation
\[
	\frac{d}{dt}P_{t} + A^{\top}P_{t} + P_{t}A -2P_{t}BB^{\top}P_{t} + F = 0
\]
with $P_{T} = F_{T}$ (see~\cite[Chapter 6]{YZ1991}). We choose the initial value from a pre-specified distribution ? and the terminal time $T=0.5$. Figure~\ref{OU_2} displays the direct comparison and $L^2$ error between the approximation $u^{\textrm{NN}}$ and the baseline $u^\ast$ of the solution to this SOC problem, which illustrates the accuracy of our DeepHJB solver.

\subsection{Metastable dynamics}
\label{sec:meta}
We consider the double well
\[
\Psi(x) = \sum_{i=1}^n\kappa_i(x_i^2 - 1)^2, \quad \kappa_i > 0.
\]
and the controlled system with
\[
	A(\mathbf{x}_t) = -\nabla\Psi, \quad B= C= I_{n\times n}, \quad \lambda =0.
\]
The initial states in this experiment are $(-1, \cdots, -1)^\top$, and the terminal state is set as $(1, \cdots, 1)^\top$. As for the cost functional, we choose
\[
F=0,\quad D = I_{n\times n}, \quad F_T = \mathrm{diag}\{\nu_1,\cdots, \nu_i, \cdots, \nu_n\}
\]
and the terminal time $T=1.0$.

\begin{figure}[!htbp]
\centering
\subfigure[Approximation $u^{\textrm{NN}}$]{
\includegraphics[width=0.3\textwidth]{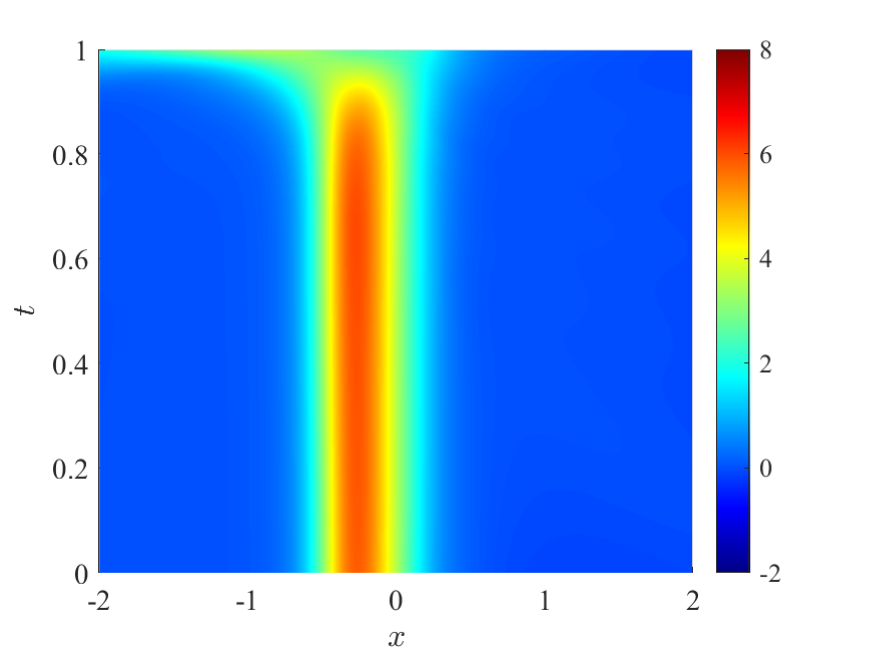}
}
\hspace{0in}
\subfigure[Baseline $u^\ast$]{
\includegraphics[width=0.3\textwidth]{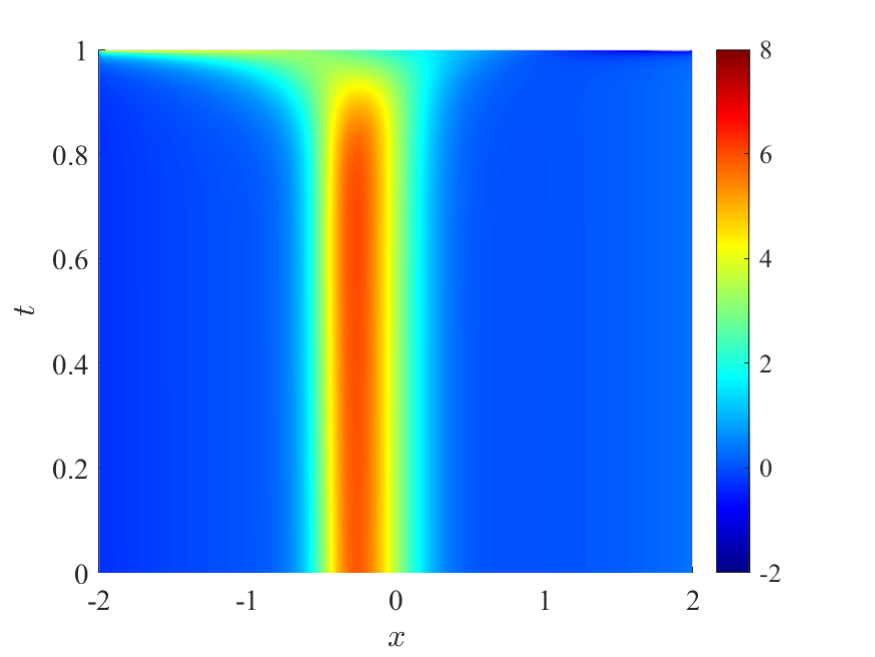}
}
\caption{Optimal control}
\label{meta1d_2}
\end{figure}

\begin{figure}[!htbp]
\centering
\includegraphics[width=0.3\textwidth]{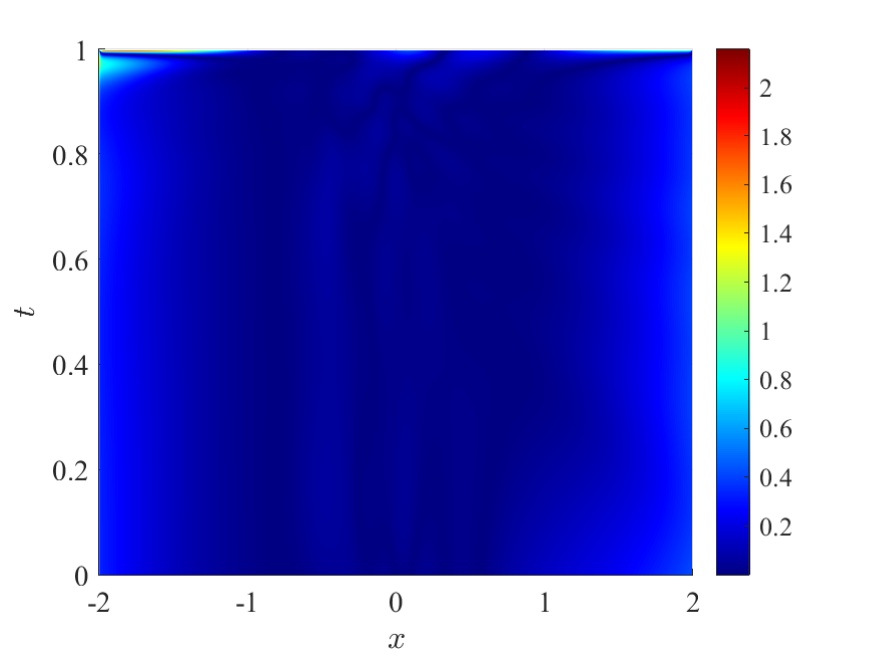}
\caption{Absolute error $|u^{\textrm{NN}} - u^\ast|$}
\label{meta1d_3}
\end{figure}

\begin{figure}[!htbp]
\centering
\includegraphics[width=0.4\textwidth]{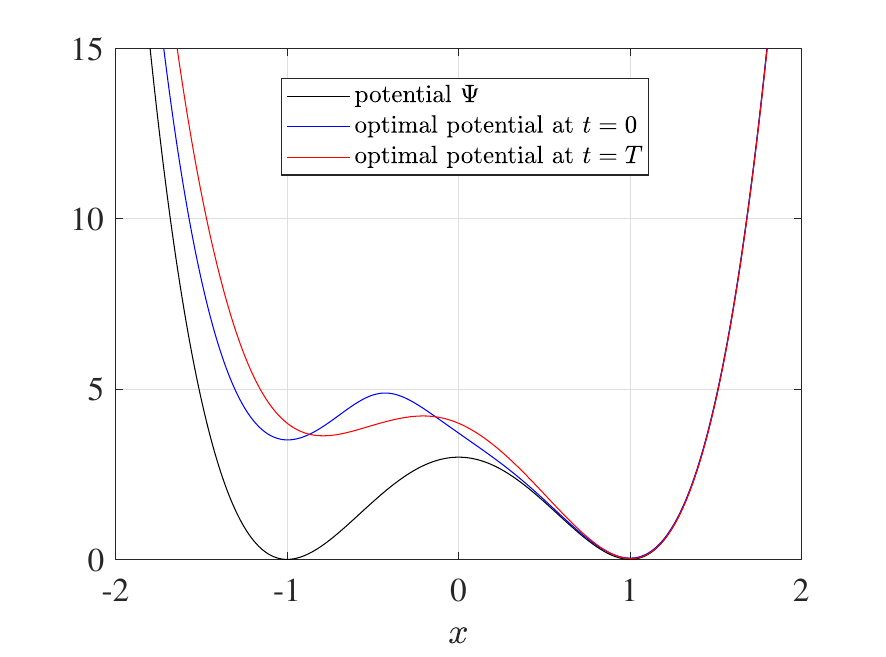}
\caption{Original potential and optimal potential.}
\label{meta1d}
\end{figure}

Firstly, we study the one-dimensional setting, choosing $\kappa = 3$, $\nu =1$. Figure~\ref{meta1d_2} displays the approximation $u^{\textrm{NN}}$ of the optimal control computed by the DeepHJB solver and the baseline $u^\ast$ obtained by a finite difference method. The absolute error between them can be seen in Figure~\ref{meta1d_3}. It is clear that the approximation is in close agreement with the baseline. Figure~\ref{meta1d} demonstrates the growth of the potential function from an original potential to the optimal potential.

\begin{figure*}[!htbp]
\centering
\includegraphics[width=0.8\textwidth]{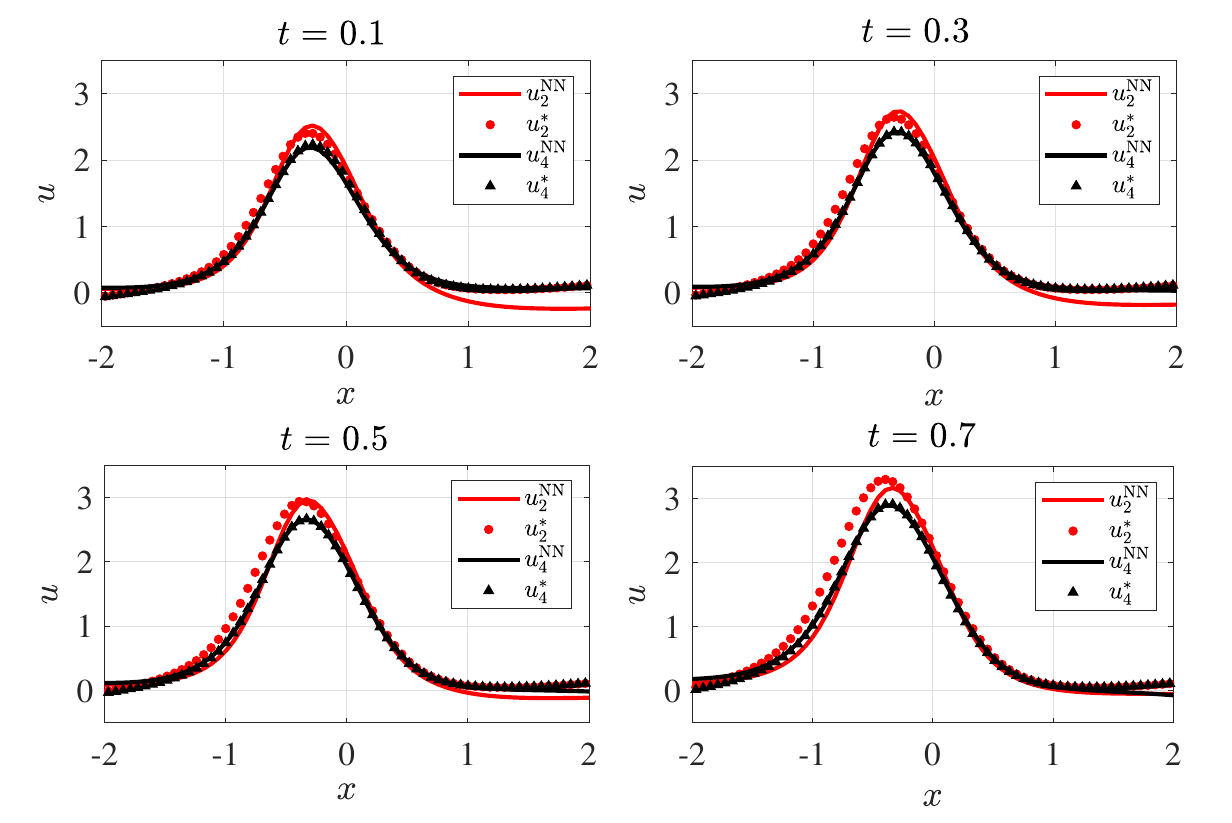}
\caption{Comparison of $u^{\textrm{NN}}$ and $u^\ast$ at different times.}
\label{meta_high}
\end{figure*}

Let us next consider the high-dimensional case, that is, $n=5$. In particular, we set $\kappa_i = 1.2$, $\nu_i = 1$ for $i \in\{1, 2, 3\}$ and $\kappa_i = 1$, $\nu_i = 1$ for $i \in\{ 4, 5\}$ .
As can be seen, Figure~\ref{meta_high} shows two components of the five dimensional approximated optimal control $u^{\textrm{NN}}$ as well as the baseline $u^\ast$, which indicates a good match and illustrates the efficacy of our DeepHJB solver for solving a high dimensional nonlinear SOC problem. 

%%%%%%%%%%%%%%%%%%%%%%%%%%%
%%%%%%%%%%%%%%%%%%%%%%%%%%%
\section{Conclusion}
\label{sec:conclusion}
In this paper, we proposed the DeepHJB solver to study the finite time horizon SOC problems for a class of dynamical systems. 
Although these numerical experiments in this work demonstrate the efficacy of the solver,  it still has plenty of room for development. 
From the viewpoint of theoretical analysis, our future research will be devoted to connecting Lyapunov analysis with the DeepHJB solver, and doing error analysis for the solver. Moreover, we will also exploit the present solver to investigate more control problems of high-dimensional nonlinear systems in practical applications. 

\appendix

\section{Proof of Theorem~\ref{verification}}
\label{sec:thmproof}
From Theorem 5.1 in Chapter 5.5 of \cite{YZ1991}, we know the fact that an admissible pair $(\mathbf{x}_t, \mathbf{u}(t))$ is optimal is equivalent to the condition that this pair satisfies the following \rm{HJB} equation
\begin{equation} \label{app:hjb}
{\small
\begin{aligned}
    	 -\partial_t q(t, \mathbf{x}_t) = &H\left(t, \mathbf{x}_t, \mathbf{u}(t), \nabla q(t, \mathbf{x}_t), \nabla^2 q(t, \mathbf{x}_t)\right) \\
	 =&\frac{1}{2}\mathrm{tr}\left[\sigma(t, \mathbf{x}_t, \mathbf{u}(t))^{ \top }\nabla^2q(t, \mathbf{x}_t) \sigma(t, \mathbf{x}_t, \mathbf{u}(t)) \right]\\
	 &+ b(t, \mathbf{x}_t, \mathbf{u}(t))\cdot \nabla q(t, \mathbf{x}_t) + \phi(t, \mathbf{x}_t, \mathbf{u}(t))
\end{aligned}
}
\end{equation}
and
\[
\mathbf{u}(t) = \mathop{\arg\min}\limits_{\mathrm{u}\in\mathcal{U}}\Big\{\frac{1}{2}\mathrm{tr}\left[\sigma^{ \top }\nabla^2q \sigma \right] + b\cdot \nabla q + \phi \Big\}.
\]
Due to the specific expression of $b$, $\sigma$ and $\phi$, we have
\begin{equation*}
\mathbf{u}(t) = \mathop{\arg\min}\limits_{u\in\mathcal{U}}\Lambda(u)
\end{equation*}
with
\begin{equation*}
\begin{aligned}
\Lambda(u):=&\frac{1}{2}\lambda^2\mathrm{tr}\left[(B(t, \mathbf{x}_t)u(t))^{ \top }\nabla^2q(t, \mathbf{x}_t) B(t, \mathbf{x}_t)u(t) \right] \\
&+ B(t, \mathbf{x}_t)u(t)\cdot \nabla q(t, \mathbf{x}_t) + \frac{1}{2}u(t)^{\top}Du(t).
\end{aligned}
\end{equation*}
Since we have $\frac{d\Lambda}{du}\big|_{u=\textbf{u}} = 0$, that is,
\begin{equation*}
\lambda^2B^{ \top }\nabla^2q B\textbf{u} + B^{ \top }\nabla q + D\textbf{u} =0,
\end{equation*}
then we have
\begin{equation} \label{app:control}
\textbf{u} = -(D+\lambda^2B^{ \top }\nabla^2q B)^{-1}B^{ \top }\nabla q.
\end{equation}
Plugging the expression of the optimal control~\eqref{app:control} into equation~\eqref{app:hjb}, we obtain the desired equation~\eqref{pathwiseValue}.

\section{Experiment configuration}
\label{sec:configuration}

We introduce the following fully connected feedforward neural network
\begin{equation*}
{\small
\begin{aligned}
z^{(1)}(x, \theta) &= W^{(1)} x + b^{(1)}, \\
\bar{z}^{(l)}(x, \theta) &= \sigma(z^{(l)}(x, \theta)), \quad l=1, 2, \cdots, L-1,\\
z^{(l+1)}(x, \theta) &= W^{(l+1)} \bar{z}^{(l)}(x, \theta) + b^{(l+1)}, \quad l=1, 2, \cdots, L-1,
\end{aligned}}
\end{equation*}
where we refer to $\sigma: \mathbb{R}\rightarrow\mathbb{R}$ as the activation function, to $L$ as the number of layers, and  to $N_0$, $N_L$, and $N_l$ as he number of neurons in the input, output, and $l$-th hidden layer, respectively. We denote by $\mathcal{A} = (N, \sigma)$, $N = (N_0, N_1, \cdots, N_L)\in\mathbb{N}^{L+1}$, the architecture of the neural network. 

The computational framework for our numerical examples is conducted by the following architecture.
\begin{itemize}
\item Controlled OU dynamics in Section \ref{sec:OUprocess}. 
\begin{itemize}
\item
For the linear terminal cost, the architecture is given by
\begin{equation*}
\begin{split}
\mathcal{A}_1 &= ((31, 32, 32, 32, 1), \tanh),\\
\mathcal{A}_2 &= ((31, 32, 32, 32, 30), \tanh),
\end{split}
\end{equation*}
%the learning rate is $0.005$, and the number of training point $M=10000$.
while for the quadratic terminal cost it is given by
 \begin{equation*}
\begin{split}
\mathcal{A}_1 &= ((16, 64, 64, 64, 1), \tanh),\\
\mathcal{A}_2 &= ((16, 64, 64, 64, 15), \tanh).
\end{split}
\end{equation*}
\end{itemize}
\item Controlled metastable dynamics in Section \ref{sec:meta}. 
\begin{itemize}
\item
For one-dimensional case, the architecture is given by
\begin{equation*}
\begin{split}
\mathcal{A}_1 &= ((2, 128, 128, 128, 128, 1), \tanh),\\
\mathcal{A}_2 &= ((2, 128, 128, 128, 128, 1), \tanh),
\end{split}
\end{equation*}
%the learning rate is $0.005$, and the number of training point $M=10000$.
while for ten-dimensional case it is given by
 \begin{equation*}
\begin{split}
\mathcal{A}_1 &= ((6, 128, 128, 128, 128, 1), \tanh),\\
\mathcal{A}_2 &= ((6, 128, 128, 128, 128, 5), \tanh).
\end{split}
\end{equation*}
\end{itemize}
\end{itemize}

The computing device that we use for our solver includes a single NVIDIA GeForce RTX 2080Ti GPU with 11GB memory. Codes will be publicly available at \url{https://github.com/zhezhejiao/DeepHJB} after being accepted.

 \bibliographystyle{elsarticle-num} 
 \bibliography{refs}

\end{document}